\documentclass[10pt]{article}

\usepackage{graphics}
\usepackage{graphicx}

\usepackage[a4paper, left=35mm,right=35mm,top=34mm,bottom=34mm]{geometry}
\usepackage[utf8]{inputenc}
\usepackage[T1]{fontenc}
\usepackage[english]{babel}

\usepackage{enumerate}
\usepackage{graphicx}
\usepackage{hyperref}
\hypersetup{
    colorlinks=true,
    linkcolor=blue,
    filecolor=magenta,      
    urlcolor=cyan,
}

\usepackage{mathtools,amsthm,amssymb,amsfonts}
\usepackage{algorithm}
\usepackage{algorithmic}
\makeatother
\theoremstyle{plain}

\theoremstyle{definition}

\theoremstyle{remark}

\usepackage{caption} 
\usepackage[tight,footnotesize]{subfigure}

\usepackage{fancyhdr}

\author{
  {\normalsize Liang S. Lai}\thanks{School of Computing and Mathematical Sciences, University of Greenwich, United Kingdom.}
	\and
	{\normalsize Choi-Hong Lai}\thanks{School of Computing and Mathematical Sciences, University of Greenwich, United Kingdom.}
	\and
  {\normalsize Abal-Kassim Cheik Ahamed}\thanks{\'Ecole Centrale Paris, France.}
	\and
  {\normalsize Fr\'ed\'eric Magoul\`es}\thanks{\'Ecole Centrale Paris, France
    (correspondence, frederic.magoules@hotmail.com).}
		}	
\title{Coupling and Simulation of Fluid-Structure Interaction Problems for Automotive Sun-roof on Graphics Processing Unit}
\date{}

\begin{document}
\maketitle
\thispagestyle{fancy}

\begin{abstract}
\noindent  In this paper, the authors propose an analysis of the frequency response function in a car compartment, subject to some fluctuating pressure distribution along the open cavity of the sun-roof at the top of a car.
Coupling of a computational fluid dynamics and of a computational acoustics code is considered to simulate the acoustic fluid-structure interaction problem.
Iterative Krylov methods and domain decomposition methods, tuned on Graphic Processing Unit (GPU), are considered to solve the acoustic problem with complex number arithmetics with double precision.
Numerical simulations illustrate the efficiency, robustness and accuracy of the proposed approaches.
\end{abstract}


\section{Introduction}

The rapid advance of computational power in recent years allows the use of large eddy simulations (LES) in many high Reynolds number applications.
The main advantage of LES over those computationally less expensive methods such as Reynolds-averaged Navier-Stokes equations (RANS) is the increased level of detail it can deliver.
While RANS methods provide ``averaged'' flow fields and over-damp high frequency fluctuations, LES is able to predict instantaneous flow characteristics and capture energy-containing eddies (i.e large turbulent scales), which are the principal contributors to sound generation in many problems.
Thus, for flows involving flow separation or acoustic prediction, LES offers significantly more accurate results over RANS approaches.
LES is also used to unravel the physics of turbulence and to compute flows of industrial relevance, when Reynolds-averaged models perform poorly or direct numerical simulation (DNS) techniques are prohibitively expensive.

In the current study, a hypothetical car configuration with an open sun-roof and a compartment forming a cavity is examined.
The car, travelling at a fixed cruising speed, experiences induced flow fluctuations due to the open sun-roof.
The pressure perturbation along the sun-roof is computed by solving the unsteady compressible Navier-Stokes equations.
For this purpose, the finite volume CFD package, PHOENICS~\cite{phoenics2014} is used.
Then these pressure fluctuations, due to the sun-roof, are extracted and analysed.
Various high order numerical schemes are compared to identify the advantages and disadvantages of each one for this application.
Furthermore, the acoustic response inside the car compartment is studied by solving the Helmholtz equation for the acoustic pressure with the ACOUFEM~\cite{acoufem2014} library.
The solution of the acoustic problem is obtained with a Krylov method, then a domain decomposition approach~\cite{magoules:journal-auth:10,magoules:journal-auth:9} which has shown strong robustness for solving acoustic problems arising from the automotive industry~\cite{magoules:journal-auth:4}.
These methods are optimized for complex number arithmetics with double precision on GPU architecture.
The interesting speed-up obtained with these methods allows a fast and accurate analysis of the acoustic phenomena within the car compartment, in an extremely short time scale.

\section{Computational fluid dynamics}

\subsection{Unsteady Navier-Stokes equations resolution}

\begin{figure}[htbp]
\centering
\includegraphics[width=7.8cm]{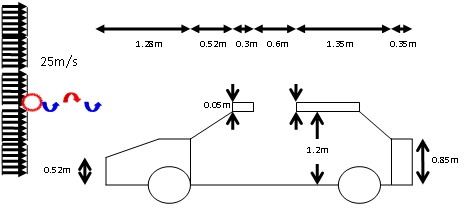}
\caption{A hypothetical car with open sun-roof}\label{Fig1}
\end{figure}

Previous experience of an open cavity with a lip shows induced oscillatory pressure fluctuations~\cite{Wang2008} caused by shear layer separation at the upstream end.
This leads to further interests in related problems such as a hypothetical car with an open sun-roof as depicted in Fig.~\ref{Fig1}.
The length of the sun-roof is equal to 0.6 $m$ and the effective depth of the opening lip (thickness of the sun-roof) is  equal to 0.05 $m$.
The free stream velocity is equal to 25 $m s^{-1}$ ($\sim$ 90 $km h^{-1}$).
To get a stronger pressure fluctuation response on top of the sun-roof, the flow is excited by an artificial sinusoidal vertical-velocity disturbance used to represent a single vortex generated by a vehicle travelling upstream of this car.
The vortex strength is given by the formula $W = W_0 \sin(2 \pi a t)$, where $W_0 = -1.2 m s^{-1}$ and $a$ is a constant parameter independent of the time.
Different frequencies of this upstream vertical-velocity disturbance are applied to generate different acoustic responses on the top of the sun-roof.

\begin{figure}[htbp]
\centering
\includegraphics[width=7.8cm]{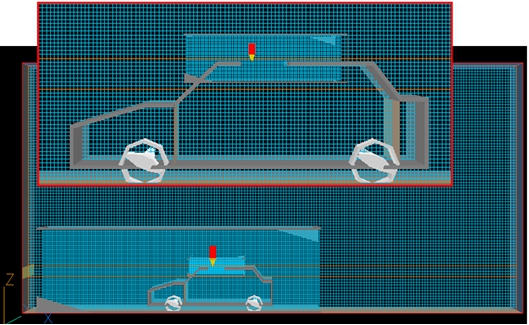}
\caption{Fine grid applied to the airflow around the car configuration and finer grid specifically focused on the top of the sun-roof opening.}\label{Fig2}
\end{figure}

In this paper, a finite volume based software package, PHOENICS~\cite{phoenics2014}, is used to compute time-accurate unsteady flow fields.
The package may be used in the computations of compressible or incompressible flows.
It uses a structured, regular Cartesian mesh and a grid refinement scheme which refines the grid size in each direction equally, as represented Fig.~\ref{Fig2}.

When a strong vortex hits the pressure outlet, it causes backflows and because of not well-converged time step solution (limited by the maximum number of iterations), it affects the vortex shedding.
That is why it is necessary to extend the domain downstream.
In the present simulation, the computational domain is taken as 17.6 $m$ by 8.8 $m$.
Four levels of grid resolution are used with a minimum grid size of 0.025 $m$.
To satisfy both the mass and momentum conservation laws, the velocity and pressure fields are solved iteratively by using the SIMPLE pressure-correction algorithm proposed by Patankar and Spalding~\cite{Patankar1980}.
In PHOENICS, standard boundary conditions are used for inflow, solid wall, and far-field boundaries.
Five different discretisation schemes have been tested in this study to provide a better understanding of their advantages and disadvantages for the present case.
In order to resolve the acoustic disturbances correctly a minimum of 20 temporal integration steps were chosen to represent each oscillation cycle at the highest frequency of interest.
The time step length, $\delta t$, chosen for the temporal integration is $10^{-3}$ $s$ resulting in a maximum resolved frequency of 50 $Hz$.

\subsection{Extraction of pressure fluctuations}

Two factors contributed to the pressure fluctuations above the sun-roof: the incoming flow over the vehicle's body and the artificial disturbance introduced upstream of the configuration.
This artificial disturbance requires a time equal to 528 $\delta t$ to reach the sun-roof's downstream for the present study.
The pressure obtained from the computational fluid dynamics (CFD) calculation is used to examine the frequency response inside the car compartment.
The pressure fluctuation along the upper surface of the car configuration and at the sun-roof opening is given by
$$P_f (x,t) = P(x,t) - \bar{P}(x,t)$$
where $P$ is the instantaneous pressure distribution along the upper surface obtained by the CFD calculation and $\bar{P}$ is the background pressure distribution along the upper surface due to the upstream velocity and the car configuration.

\subsection{Numerical schemes} 

\begin{figure}[htbp]
\centering
\includegraphics[width=7.8cm]{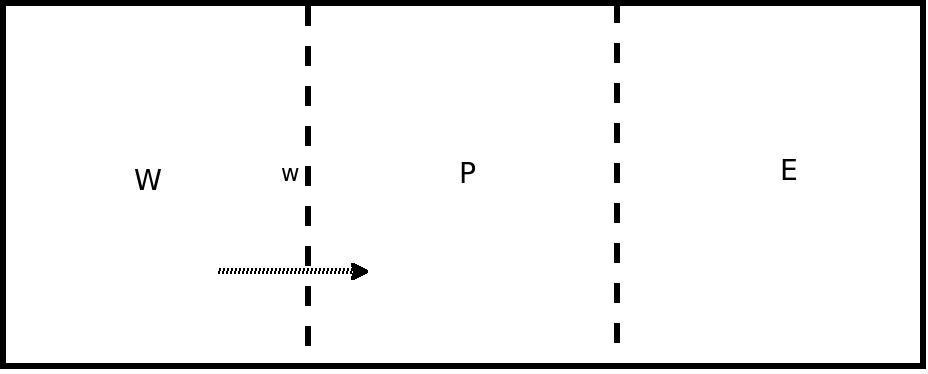}
\caption{Centered finite volume method}\label{FigX}
\end{figure}

In all cell centered finite volume methods, as illustrated in Fig.~\ref{FigX}, values of the variable $\phi$ are known at the cell centres $W$, $P$ and $E$; but the values of $\phi$ at face $w$ are not known and may be calculated by using a number of numerical schemes.

The numerical scheme influences the balance equations for both cell $W$ and cell $P$.
To obtain a fairly good solution one can choose $\phi_w = \phi_W$ when the flow goes from $W$ to $P$ or $\phi_w = \phi_P$ when the flow travels from $P$ to $W$.
This kind of scheme is used as the default numerical scheme, together with other schemes, in PHOENICS.
In this paper, five different numerical schemes, three linear and two non-linear schemes, as listed below, are used.
Each of them has a different approach to calculate the cell face value $\phi_w$ .
\begin{itemize}
\item Upwind-differencing scheme (UDS): \\
$$\phi_w = \phi_W$$
\item Central-differencing scheme (CDS): \\
$$\phi_w = \frac{\phi_W + \phi_P}{2}$$
\item Quadratic upwind scheme (QUICK): \\
$$\phi_w = \frac{3}{8}\phi_P + \frac{3}{4}\phi_W - \frac{1}{8}\phi_{WW}$$
\item SMART Bounded QUICK: \\
$$\phi_w = \phi_W + 0.5 B (\phi_W-\phi_{WW})$$
\item Harmonic QUICK (HQUICK): \\
$$\left \{ \begin{array} {ll}
\phi_w = \phi_W & \textrm{ if  $r \leq 0 $}, \\ [0.2cm]
\phi_w = \phi_W + \frac{2(\phi_P-\phi_{W})(\phi_W-\phi_{WW})}{\phi_P + 2\phi_{W} - 3\phi_{WW}} & \textrm{if $r > 0 $}. \\
\end{array} \right.$$
\end{itemize}
where
$$B = \max \left( 0 , \min \left( 2r, \frac{3r+1}{4} , 4 \right) \right), r= \displaystyle{\frac{\phi_P - \phi_W}{\phi_W - \phi_{WW}}}$$
and $\phi_{WW}$ is the cell-centred value of $\phi$ at the cell upstream the cell W.

\begin{figure}[htbp]
\centering
\includegraphics[width=7.8cm]{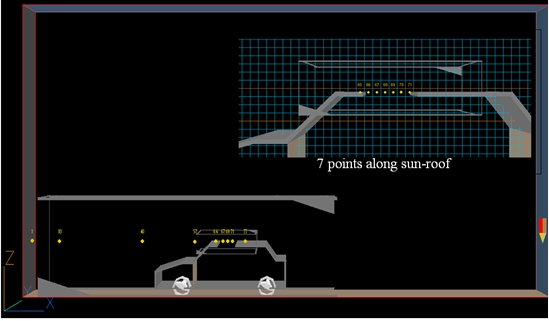}
\caption{Nine observation points in the computational domain. On top: seven points along the sun-roof}\label{Fig3}
\end{figure}

Fig.~\ref{Fig3} shows nine observation points marked with their node numbers, along the line $y = 1.6$ $m$, i.e. just one-cell above the sun-roof on the same streamline where the upstream artificial disturbance is introduced.
The time history of the pressure fluctuations at these observation points are shown in column~1 of Table~\ref{Table1}.
A zoom-in to the neighbourhood of the sun-roof using seven other observation points for pressure fluctuations are shown in column~2 of Table~\ref{Table1}.
It can easily be seen that first order accurate Hybrid/Upwind scheme and third order accurate HQUICK scheme are too dissipative and therefore are not suitable for this type of example.
As a result, the magnitude of pressure fluctuations observed on top of the sun-roof is very small.
CDS even failed to converge because the cell Peclet number is not guaranteed to be less than two.
However, SMART and QUICK scheme show more interesting results.
The pressure fluctuations on top of the sun-roof gradually grow in magnitude in a sinusoidal form.
The fluctuations obtained by using QUICK scheme lead to a more stable and regular sinusoidal shape.
\begin{table}[htbp]
\begin{center}
\caption{Comparison of pressure fluctuations using different numerical schemes. From top to bottom Hybrid \/ Upwind, HQUICK, SMART, and QUICK schemes.}\label{Table1}
\begin{tabular}{cc}
Pressure fluctuation & Pressure fluctuation  \\
at nine              &  at seven points \\
observation points   &  on top of sun-roof \\
\includegraphics[scale=0.79]{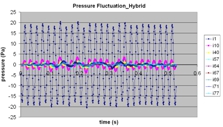} & \includegraphics[scale=0.79]{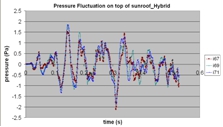} \\
\includegraphics[scale=0.79]{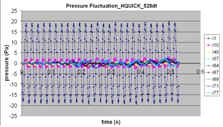} & \includegraphics[scale=0.79]{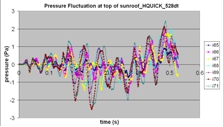} \\
\includegraphics[scale=0.79]{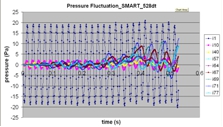} & \includegraphics[scale=0.79]{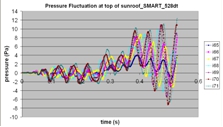} \\
\includegraphics[scale=0.79]{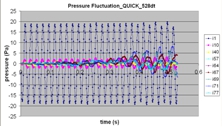} & \includegraphics[scale=0.79]{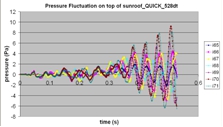} \\
\end{tabular}
\end{center}
\end{table}

\begin{figure}[htbp]
\centering
\includegraphics[width=7.8cm]{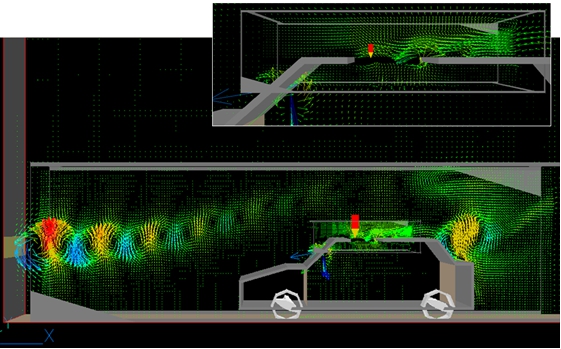}
\caption{A snapshot of $z$ component velocity disturbance at $t = 5\,s$. At the bottom: the zoom-in image on top of the sun-roof}\label{Fig4}
\end{figure}

On Fig.~\ref{Fig4} a snapshot of vertical velocity disturbance is shown.
It points out that the amplitude of aerodynamic disturbances becomes weaker and weaker because of the numerical scheme dissipation.
However, a clear vortex shearing on top the sun-roof can still be observed.
At this stage, QUICK seems to be the best high order scheme to be applied for this application.

\section{Acoustic analysis}

\subsection{Frequency response function}

\begin{figure}[htbp]
\centering
\includegraphics[width=7.8cm]{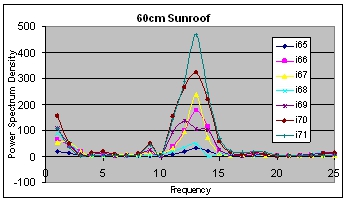}
\caption{Power spectrum density of the time history via a 512-point FFT}\label{Fig5}
\end{figure}

Frequency components of the pressure fluctuations are then examined by computing the acoustic power spectrum for the temporal signals monitored at all seven points on the sun-roof via a 512-point Fast Fourier Transform (FFT).
The spectrum, reproduced on Fig.~\ref{Fig5}, shows that the dominant frequency at all observation points on the sun-roof roughly occurs at 13 $Hz$.
The validity of the results for the dominant frequency is checked using a Helmholtz resonator with a geometrically similar cavity.
The resonant frequency for a typical Helmholtz resonator can be approximated by the formula, 
$$
\displaystyle{ f = \frac{c}{2 \pi}\sqrt{\frac{A}{l_{eff}V}}},
$$
where $l_{eff} = l + l_{cor} + \eta \, r$ denotes the effective length of the air at the opening, $l$ is the geometric neck length (i.e., 0.05 $m$, on Fig.~\ref{Fig1}).
$l_{cor}$ is the end correction on the neck length, which can be expressed by a product of $r$, the radius of the neck, and $\eta$, an empirical coefficient which significantly depends on the geometrical configuration.
$A$ is the cross sectional area of the neck and $V$ represents the volume of the inside cavity.
Although it is an idealized formula that completely neglects the shear layer, it gives an approximate indication for the oscillation frequency of the cavity.
An approximate value of the dominant resonant frequency with $\eta = 16.9$ is around 10.5 $Hz$.
It must be pointed out that this is not a strict comparison as the coefficient $\eta$ is currently unavailable for the cavity formed by the considered car compartment.
However, even so, this crude comparison shows that the dominant frequency value obtained through the unsteady computation is a physically acceptable approximation.

To study the acoustic response along the sun-roof, different frequencies of the upstream disturbance are applied.
Numerical tests as a function of input disturbance are performed to verify the hypothesis that the lower the frequency of disturbance, the lower the frequency of acoustic response obtained.
In this study, 25 $Hz$ and 10 $Hz$ disturbance frequencies are compared with the maximum resolvable frequency of 50 $Hz$ and are reported in Table~\ref{Table2}.
The power spectra shows that for an incoming disturbance at a frequency higher than 25 $Hz$ the dominant mode of the generated noise due to the sun-roof roughly occurs at 13 $Hz$ which is the resonant frequency.
On the other hand a lower frequency incoming disturbance, say at 10 $Hz$, seems to excite a half harmonic at around 6 $Hz$ while maintaining the fondamental harmonic (13 $Hz$) at a weaker strength.
\begin{table}[htbp]
\begin{center}
\caption{Comparison of peak frequencies obtained via a 512-point FFT due to different incoming disturbances. From top to bottom at 50, 25 and 10 $Hz$}\label{Table2}
\begin{tabular}{c}
\includegraphics[scale=0.7]{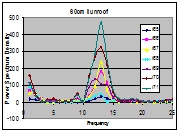} 
\includegraphics[scale=0.7]{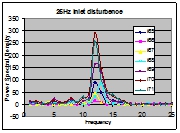} 
\includegraphics[scale=0.7]{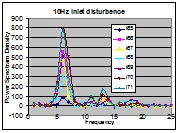}
\end{tabular}
\end{center}
\end{table}

\subsection{Wave equation inside the car compartment}

In Lighthill's acoustic analogy~\cite{Lighthill1952}, Lighthill essentially recasts the exact equations of fluid motion (Navier Stokes equations and continuity equation) in the form of an inhomogeneous wave equation suitable to be applied in the far-field acoustics, therefore making an acoustic analogy with fluid mechanics.
The conservative form of the continuity equation~(\ref{eq1}) and momentum equation~(\ref{eq2}) for a compressible fluid in a Cartesian coordinate system, without body forces are
\begin{eqnarray}
\frac{\partial \rho}{\partial t} + \frac{\rho u_i}{\partial y_j} & = & 0 ,\label{eq1} \\
\frac{\partial (\rho \, u_i)}{\partial t} + \frac{\rho \, u_i \, u_j + p_{ij}}{\partial y_j} & = & 0 , \label{eq2}
\end{eqnarray}
where $p_{ij} = p\delta_{ij} - \tau_{ij}$. $\rho$ is the density, $\delta_{ij}$ is the Kronecker symbol (i.e., $\delta_{ij} = $ if $i = j$ and $\delta_{ij} = 0$ otherwise), $u_i$, $u_j$ are the velocity components, $p_{ij}$ is the stress tensor and $p$ is the static pressure.
If external sources are ignored, the famous Lighthill's wave equation can be written as 
\begin{equation}
\frac{\partial^{2} \rho^{'}}{\partial^{2} t} - c_0^{2} \, \nabla^{2} \, \rho^{'} = \frac{\partial^{2} \, T_{ij} }{\partial y_j \, \partial y_j}, \label{Lighthill}
\end{equation}
where $\rho^{'}$ is density perturbation (defined as the deviation from the quiescent reference density) and $c_0$ is the speed of sound in the fluid at rest.
The Lighthill stress tensor, $T_{ij}$ is defined as
\begin{equation}
T_{ij} = \rho \, u_i \, u_j - \tau_{ij} + (p^{'} - c_0^{2}\rho^{'})\delta_{ij}, \label{Tensorlight}
\end{equation}
where $p^{'}$ is the pressure perturbation and $\tau_{ij}$ is the viscous stress term.
Each of these acoustic source terms may play a significant role in noise generation depending on the considered flow conditions.
It is however generally accepted, that the term $\tau_{ij}$ has a lower impact, by several orders of magnitude, on noise generation than the other terms and can consequently be neglected in most situations.
Note that the perturbations ($p^{'}$,$\rho^{'}$) are defined as the deviations between the total flow variables ($p$,$\rho$) and the quiescent reference state ($p_0$,$\rho_0$) during the derivation of equation~(\ref{Lighthill}), i.e.,
\begin{eqnarray}
p^{'} & = & p - p_0, \\
\rho^{'} & = & \rho - \rho_0.
\end{eqnarray}
The source term on the right hand side of equation~(\ref{Lighthill}) consists in a detailed flow motion around the acoustic source region (near-field).
In this particular case, as the ratio between the width of the opening of the cavity (i.e. sun-roof, in this case) and the depth of the cavity (i.e. height inside car compartment) is relatively high, the aerodynamic motion inside the car compartment can be effectively neglected.
In other words, Lighthill's equation can be rewritten as a homogeneous wave equation expressed as
\begin{equation}
\frac{\partial^{2} p^{'}}{\partial^{2} t} - c^{2} \, \nabla^{2} \, p^{'} = 0, \label{pressure}
\end{equation}
where $p^{'}$ is pressure fluctuation.
In order to transform the wave equation from time domain to frequency domain, we need to integrate equation~(\ref{pressure}) with respect to time by using a fast Fourier transform (FFT),
\begin{equation}
\int^{\infty}_{-\infty} \, \frac{\partial^{2} p^{'}}{\partial^{2} t} - c^{2} \int^{\infty}_{-\infty} \, \nabla^{2} \, p^{'} = 0, \label{pressurefft}
\end{equation}
and one can get in a simplified form the Helmholtz equation~(\ref{Helmoltz}),
\begin{equation}
- \omega^2 \, \psi - c^2 \, \nabla^2 \psi = 0 \mbox{, where} \, \psi = \int^{\infty}_{-\infty} p^{'} e^{i \, \omega \, t}. \label{Helmoltz}
\end{equation}
It is assumed that the flow inside the car compartment is negligible, thus the acoustic propagation can be described by the Helmholtz equation.
The analysis of the sound distribution in the car compartment for the dominant frequency ($f =$ 13 $Hz$) due to an incoming disturbance of 50 $Hz$ is carried out.
The power spectrum density along the sun-roof is used as a Dirichlet boundary condition for the Helmholtz equation, which outputs the acoustic pressure distribution inside the car compartment.
A unstructured mesh is considered for the car compartement and stabilized finite elements are considered for the discretization~\cite{magoules:journal-auth:7}.

\begin{figure}[htbp]
\centering
\includegraphics[width=7cm]{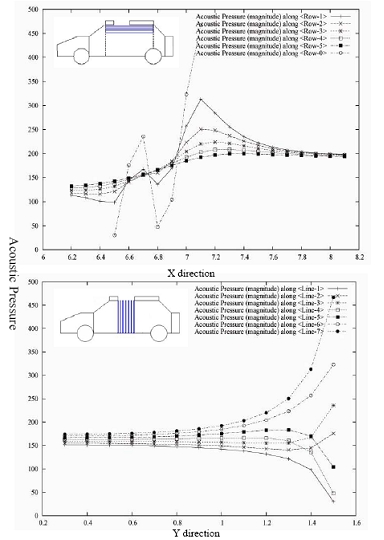}
\caption{Acoustic pressure inside the car component along several horizontal and vertical lines}\label{Fig6}
\end{figure}

Fig.~\ref{Fig6} represents the acoustic pressure distribution along several horizontal and vertical lines inside the car compartment.
It shows that the highest acoustic pressure is experienced at 7.1 $m$ in the $x$ direction.
On the other hand along the horizontal line just below the sun-roof the pressure shows an oscillatory behaviour resulting from the pressure fluctuation above the sun-roof.
This oscillatory behaviour is gradually weakened as one moves deeper into the car compartment.
The acoustic pressure distribution along all vertical lines seems to show that oscillatory effects disappear deep inside the car compartment.
These observations show that the obtained solution is reliable.
The acoustic pressure tends to be more stable at the bottom of the car compartment.

\section{GPU computations}

With the apparition of Graphics Processing Unit (GPU) in 2000s, parallel computations have been greatly enhanced.
Indeed, performances of Central Processing Unit (CPU) and GPUs are significanlty different, due to inherently different architectures.
GPU architecture exploits parallelism by having many floating points processors exploiting large amount of data in parallel.
In order to perform fast computations, special memory hierarchy allowing each processor to optimally access the requiered data is mandatory.
The memory is thus a key feature of GPU architecture~\cite{ref:cheikahamed:2012:inproceedings-1:DARG:2013}~\cite{ref:cheikahamed:2012:inproceedings-1:ZZ:2013}.
In summary, a CPU is constantly accessing the random-access memory (RAM), implying a low latency at the detriment of its raw throughput~\cite{ref:cheikahamed:2012:inproceedings-1:DARG:2013:b}.
Opposite, a GPU has four main types of memory: (i) the global memory that ensures the interaction with the host (CPU)--this is the slowest memory; (ii) the constant memory that provides interaction with the host--this memory is generally cached for fast access; (iii) the shared memory that is accessible by any thread\footnote{A thread is the smallest sequence of programmed instructions that can be managed independently by an operating system.} of the block\footnote{Threads are grouped into blocks and executed in parallel simultaneously. Threads are grouped into blocks and blocks are grouped together in a grid.} from which it was created--this memory is much faster than the global memory; (iv) the local memory that is specific to each compute unit.

Here, the major problem related to the implementation on GPU is that the acoustic problem involves complex number arithmetics.
Indeed, GPU have initially been developped for integer artihmetics, and using floating point operations with real number reduce significantly their performance.
This is even more problematic when using double precision floating point operations which imply a drop in performance.
Since the acoustic problem involves complex number artihmetics with double precision floating point operations, the expected performance are extremely low.

In addition, the acoustic problem to solve reduces to the solution of a large size sparse linear system.
Iterative Krylov methods are well suited for this type of problem; such methods require linear algebra operations such as sparse matrix-vector product which is definitely the most time consuming operation.
In this study, the acoustic matrix is stored in compressed sparse row (CSR) format in order to optimize the memory storage and to make advantage of sparse structure for memory access within the GPU.
In order to program the GPU, we use the CUDA language.

As mentionned earlier, CUDA was initially devoted for real numbers arithmetics.
However, given that a complex number is a set of two real numbers (real part, imaginary part), implementation is possible by defining a structure containing two real numbers.
The library provided by Nvidia proposes a structure called \verb+cuComplex+, but for performance issues, we define our own complex class template structure \verb+complex<T>+, which redefines all the operations given by \verb+std::complex+.
In order to take the best advantage of GPU architecture, the basic linear operation kernel need also to be rewritten~\cite{journals:KnibbeOV11},~\cite{GPU:BG:2009},~\cite{GPU:BG:2008}.
Former analysis~\cite{magoules:proceedings-auth:44} performed on real number artihmetics with double precision allowed to obtain, with a suitable rewritting of the CUDA kernel, excellent speed-up performance for basic linear operations and for Krylov methods~\cite{magoules:proceedings-auth:45,magoules:proceedings-auth:48}.
As a consequence, we decide to extend this analysis, to the present acoustic problem, i.e., to develop efficient iterative Krylov methods for solving linear systems with complex number arithmetics.
We have thus developped a preconditionned bi-conjugate gradient stabilized method (P-Bi-CGSTAB), a preconditionned P-BiCGSTAB parametered (l) and a preconditionned transpose-free quasi-minimal residual method (P-tfQMR)~\cite{saad_iterative_2003}, all of them with optimized CUDA kernel and with dynamic auto-tuning on GPU.
As indicated in~\cite{magoules:proceedings-auth:48} for real number arithmetics, our template implementation outperforms Cusp~\cite{GPU:CUSP:2010}, CUBLAS~\cite{GPU:CUBLAS}, CUSPARSE~\cite{GPU:CUSPARSE4.0:2011}, but behavior for complex number arithmetics with double precission is still a challenge.

Table~\ref{Table2} shows the speed-up obtained with our GPU implementation compared to the CPU implementation: the first column presents the size of the mesh, the second column gives the number of iterations, the third and the fourth columns represents respectively the CPU and GPU time in seconds, and the last column collects the speed-up.
The residual threshold is fixed to $10^{-9}$.
\begin{table}[htbp]
\setlength\columnsep{0.1pt}
\begin{center}
\caption{Speed-up of the acoustic solver (complex number arithmetics with double precision)}\label{Table2}
\begin{tabular}{lcccc}
\hline
size (h) & \#iter & CPU time (s) & GPU time (s) & speed-up \\
\hline
\multicolumn{3}{l}{\emph{P-BiCGSTAB}} & & \\
0.133425 & 21 & 0.01 & 0.030 & \textbf{0.33} \\
0.066604 & 53 &  0.24 & 0.106 & \textbf{2.26} \\
0.033289 & 94 &  4.01 & 0.703 & \textbf{5.71} \\
0.016643 & 183 & 85.70 & 9.209 & \textbf{9.31} \\
\hline
\multicolumn{3}{l}{\emph{P-BiCGSTAB(8)}} & & \\
0.133425 & 6 & 0.03 & 0.110 & \textbf{0.27} \\
0.066604 & 12 & 0.52 & 0.286 & \textbf{1.82} \\
0.033289 & 31 & 12.47 & 2.162 & \textbf{5.77} \\
0.016643 & 70 & 266.26 & 30.100 & \textbf{8.85} \\
\hline
\multicolumn{3}{l}{\emph{P-QMR}} & & \\
0.133425 & 24 & 0.02 & 0.040 & \textbf{0.50} \\
0.066604 & 52 & 0.27 & 0.113 & \textbf{2.40} \\
0.033289 & 99 &  4.71 & 0.755 & \textbf{6.24} \\
0.016643 & 214 &  102.17 & 10.786 & \textbf{9.47} \\
\hline
\end{tabular}
\end{center}
\end{table}
When we refine the mesh, i.e, that the size of the problem increases, the results become more accurate since more details in the car compartment are taken into account.
As we can see in Table~\ref{Table2}, for all our implementations the speed-up increases when the size of the problem increases.
Nevertheless, when the size of the problem becomes too large for GPU memory, which is often very limited on most GPUs, other methods must be considered.
Domain decomposition methods~\cite{Smith:1996:DPM},~\cite{quarteroni_domain_1999},~\cite{toselli_domain_2004},~\cite{magoules:journal-auth:21} based on iterative algorithms are an alternative. 
Such methods have encountered strong success for the solution of coercive elliptic problems~\cite{NME:NME1620320604},~\cite{cai_overlapping_1998} and are easy to implement on parallel computers.  
Using absorbing boundary transmission conditions on the interface between the subdomains~\cite{magoules:journal-auth:16} is a key point to obtain a fast convergence of the  domain decomposition algorithm, such as the Schwarz algorithm~\cite{magoulesf_contrib_3:Lions:1988:SAM},~\cite{magoulesf_contrib_3:Lions:1989:SAM},~\cite{magoulesf_contrib_3:Lions:1990:SAM}.
First works presented in references~\cite{Despres:1993:DDM},~\cite{magoulesf_contrib_3:Despres:1992:DHM},~\cite{Benamou:1997:DDM} consider Robin type absorbing boundary transmission conditions on the interface, which have then been optimized with a continuous approach in~\cite{magoulesf_contrib_3:Chevalier:1998:SMO}, \cite{magoulesf_contrib_3:Gander:2000:OSM}, \cite{magoules:journal-auth:24}, \cite{magoules:journal-auth:23}, \cite{magoules:journal-auth:18}, \cite{magoules:journal-auth:14}.
Further works have considered a discrete optimization of the interface conditions as first introduced in~\cite{magoules:proceedings-auth:6}, and then in~\cite{magoules:journal-auth:29}, \cite{magoules:journal-auth:20}, \cite{magoules:journal-auth:17}, \cite{magoules:journal-auth:12}, \cite{magoules:journal-auth:30}.

GPU implementation of domain decomposition method consists to allocate each subdomain to one single processor (CPU) and to solve the linear system on GPU at each iteration, as first proposed for the FETI method in~\cite{Papadrakakis:2011} and for the optimized Schwarz method in~\cite{magoules:proceedings-auth:50}.
Here the same methodology is considered, but for complex number arithmetics with double precision, and the optimized Schwarz method with two sided interface conditions is considered~\cite{magoules:journal-auth:28}.
The acoustic problem solved with the optimized Schwarz method with two sided interface conditions gives a speed-up of {\bf 7.03} on the problem with a mesh size equal to 0.008321; which is an excelent result according to the Amdahl's Law~\cite{amdahl} applied to domain decomposition method.

\section{Conclusion}

Sound distribution inside a car compartment due to incoming disturbance over the sun-roof is here studied through a coupled acoustic fluid-structure analysis.
Input disturbance has a strong influence on the particular harmonic excited.
The car compartment responds at a particular frequency of oscillation of 13 $Hz$, and lower frequency input excites a half harmonic at roughly 6 $Hz$.
Higher order schemes are necessary to extract such pressure fluctuations.

Due to the size of the acoustic problem, iterative Krylov methods and a domain decomposition method are considered.
Implementation on graphic processing unit with CUDA language, leads to excelent speed-up as soon as the CUDA kernel are rewritten and optimized for complex number arithmetics with double precision.

\section*{Acknowledgements}
The authors acknowledge the CUDA Research Center at Ecole Centrale Paris (France) for its support and for providing the computing facilities.

\bibliography{bib/GPUComputing,bib/Vibroacousts,bib/MAGOULES-PROCEEDINGS1,bib/MAGOULES-PROCEEDINGS2,bib/MAGOULES-JOURNAL1,bib/MAGOULES-CHAPTER1,bib/magoulesf_contrib_3}
\bibliographystyle{abbrv}

\end{document}